\documentclass[12pt]{amsart}
\usepackage{amsfonts}

\usepackage{graphics}
\usepackage{amsmath}
\usepackage{amscd}

\newtheorem{proposition}{Proposition}[section]
\newtheorem{prop}{Proposition}[section]
\newtheorem{theorem}[proposition]{Theorem}
\newtheorem{lemma}[proposition]{Lemma}
\newtheorem{corollary}[proposition]{Corollary}

\newtheorem{definition}[proposition]{Definition}

\newtheorem{Ack}[proposition]{Acknowledgments}
\newtheorem{ques}[proposition]{Question}

\begin{document}
\author{Emmanuel Breuillard}
\date{January 5, 2006}
\title{On uniform exponential growth for solvable groups}
\address{Emmanuel Breuillard, UFR de Math\'{e}matiques, Universit\'{e} de Lille,
59655 Villeneuve d'Ascq, France}
\email{emmanuel.breuillard@math.univ-lille1.fr}

\begin{abstract}
Using a theorem of J. Groves we give a ping-pong proof of Osin's uniform
exponential growth for solvable groups. We discuss slow exponential growth
and show that this phenomenon disappears as one passes to a finite index
subgroup.
\end{abstract}
\maketitle

\section{Introduction}

The main purpose of the present note is to give an alternative proof as well
as a strengthening of the fact, proved by Alperin \cite{Alp} (polycyclic
case) and Osin \cite{Os, Os2} that finitely generated solvable groups of
exponential growth have uniform exponential growth. Our approach is quite
different from the one taken up in those works and relies on a direct
ping-pong argument.

Let $\Gamma $ be a group generated by a finite subset $\Sigma $. Assume that 
$\Sigma $ is symmetric (i.e. $s\in \Sigma \Rightarrow s^{-1}\in \Sigma $),
contains the identity $e$, and let $\mathcal{G}=\mathcal{G}(\Gamma ,\Sigma )$
be the associated Cayley graph. The set $\Sigma ^{n}$ is the set of all
products of at most $n$ elements from $\Sigma ,$ i.e. the ball of radius $n$
centered at the identity in $\mathcal{G}$ for the word metric. Let $\mathcal{%
C}$ be the set of all such finite generating subsets $\Sigma $. We introduce
the following definition:

\begin{definition}
Two elements in a group are said to be \textbf{positively independent} if
they freely generate a free semigroup. The \textbf{diameter of positive
independence }of a Cayley graph $\mathcal{G}(\Gamma ,\Sigma )$ is the
quantity $d^{+}(\Sigma )=\inf \{n\in \Bbb{N}$, $\Sigma ^{n}$ contains two
positively independent elements$\}.$ Similarly, the diameter of positive
independence of the group $\Gamma $ is defined by $d_{\Gamma }^{+}=\sup
\{d^{+}(\Sigma ),$ $\Sigma \in \mathcal{C}\}.$
\end{definition}

The next definition is more standard:

\begin{definition}
Assume that $\Gamma $ is finitely generated. For $\Sigma $ in $\mathcal{C}$
we can define the \textbf{algebraic entropy} of the pair $(\Gamma ,\Sigma )$
to be the quantity $S_{\Gamma }(\Sigma )=\lim \frac{1}{n}\log (\#\Sigma
^{n}).$ Similarly, the algebraic entropy of $\Gamma $ is defined by $%
S_{\Gamma }=\inf_{\Sigma \in \mathcal{C}}S_{\Gamma }(\Sigma ).$
\end{definition}

It is easy to see that $S_{\Gamma }(\Sigma )$ is either positive for all $%
\Sigma $ in $\mathcal{C}$ or $0$ for all $\Sigma $ simultaneously.
Accordingly, the group $\Gamma $ is said to have exponential or
sub-exponential growth. If $S_{\Gamma }>0$, then $\Gamma $ is said to have 
\textit{uniform exponential growth}. If two elements generate a free
subsemigroup, there are exactly $2^{n}$ elements that can be written as
positive words of length $n$ in these two elements, hence the latter
quantities are related by the following inequality: 
\begin{equation*}
S_{\Gamma }\geq \frac{\log 2}{d_{\Gamma }^{+}}
\end{equation*}

In particular, if $\Gamma $ has a finite positive independence diameter,
i.e. $d_{\Gamma }^{+}<+\infty ,$ then $\Gamma $ has uniform exponential
growth, i.e. $S_{\Gamma }>0$. No converse is valid in general. As was proved
by Osin in \cite{Os3}, the free Burnside groups of large odd exponent are
uniformly non-amenable and in particular have uniform exponential growth.
However, these groups obviously have an infinite diameter of positive
independence.

In their seminal papers \cite{Mil, W} J. Milnor and J. Wolf proved that
finitely generated groups that are finite extensions of a nilpotent group
(i.e. virtually nilpotent groups) have polynomial growth (hence $S_{\Gamma
}=0$), while finitely generated solvable but non virtually nilpotent groups
have exponential growth. Refining their methods, J. Rosenblatt showed
subsequently in \cite{Ros} that any finitely generated solvable group
contains a free semigroup on two generators unless it is virtually
nilpotent, and C. Chou \cite{Chou} extended this dichotomy to the class of
elementary amenable groups.

In the eighties, M. Gromov asked whether a group with exponential growth
must have uniform exponential growth. Recently, J. Wilson \cite{Wil}
answered the question negatively by constructing several examples of
finitely generated subgroups of the automorphism group of a rooted tree such
that $S_{\Gamma }=0$ although they contain a free subgroup and hence have
exponential growth.

However, D. Osin \cite{Os, Os2} (and independently R. Alperin \cite{Alp} in
the polycyclic case) proved that $S_{\Gamma }>0$ for finitely generated
solvable or even elementary amenable groups that are not virtually
nilpotent. The class of elementary amenable groups is the smallest class of
groups containing both abelian groups and finite groups and that is stable
under subgroups, quotients, extensions and direct limit (see Chou \cite{Chou}%
). Obviously any solvable group is elementary amenable. Although it is not
explicitely stated in Osin's paper \cite{Os2}, reading between the lines of
his proof one can fairly simply derive (see Section 6) the following
stronger statement.

\begin{theorem}
\label{first}Let $\Gamma $ be a finitely generated non virtually nilpotent
elementary amenable group. Then $\Gamma $ has finite positive independence
diameter, i.e. $d_{\Gamma }^{+}<+\infty $.
\end{theorem}

The idea of proving the finiteness of $d_{\Gamma }^{+}$ in order to obtain
uniform exponential growth has been used in many places in the past such as
the work of Eskin-Mozes-Oh on non virtually solvable linear groups of
characteristic zero \cite{EMO}, or for hyperbolic groups in Gromov's
original monography \cite{Grom}. Note finally that the existence of groups
of intermediate growth shows that the above theorem cannot be extended to
all amenable groups.

In this paper we want to address the following question: how large can $%
d_{\Gamma }^{+}$ be? The following theorem and corollary give an upper bound
on $d_{\Gamma }^{+}$ for finitely generated solvable groups.

\smallskip

\begin{theorem}
\label{main}Let $\Gamma $ be a finitely generated solvable group. Then one
(and only one) of the following is true:

$(i)$ $\Gamma $ is virtually nilpotent (i.e. it contains a nilpotent
subgroup of finite index).

$(ii)$ $\Gamma $ has a finite index subgroup $\Gamma _{0}$ such that $%
d_{\Gamma _{1}}^{+}\leq 3$ for any finite index subgroup $\Gamma _{1}\leq
\Gamma _{0}.$
\end{theorem}

It is easy to check that a subgroup $\Gamma _{0}$ of index $n$ in a finitely
generated group $\Gamma ,$ satisfies $d_{\Gamma }^{+}\leq (2n+1)\cdot
d_{\Gamma_{0}}^{+}$. Hence we obtain:

\begin{corollary}
\label{cor} Let $\Gamma $ be a finitely generated solvable group which is
not virtually nilpotent. Then there is a number $C(\Gamma )$ such that $%
d_{\Gamma ^{\prime }}^{+}\leq C(\Gamma )$ for every finite index subgroup $%
\Gamma ^{\prime }$ of $\Gamma $.
\end{corollary}

Hence within a given commensurability class of finitely generated non
virtually nilpotent solvable groups, one can always find a group $\Gamma $
with universally bounded growth, i.e. $S_{\Gamma }\geq \frac{\log 2}{3}.$ In
particular, slow exponential growth is a phenomenon that disappears
completely as one passes to a suitable finite index subgroup.

The method used here to prove Theorem \ref{main} is based on a direct
ping-pong argument and hence differs radically from those of the above
mentioned previous works \cite{Mil, W, Ros, Os, Os2}. Thanks to the
following results that can be derived from the work of J. Groves, it is
enough to prove Theorem \ref{main} for metabelian groups (i.e. extensions of
two abelian groups) and even for subgroups of affine transformations of a $K$%
-line. This on the other hand is fairly simple as explained in Section 3.

\begin{theorem}
\label{JNVN}(\cite{Gro}) Let $\Gamma $ be a non virtually nilpotent finitely
generated solvable group all of whose proper quotients are virtually
nilpotent. Then we have:

$(i)$ $\Gamma $ is virtually metabelian.

$(ii)$ If $\Gamma $ is metabelian, then it embeds in the group of affine
transformations of the $K$-line, for some field $K$.
\end{theorem}

For the convenience of the reader, we will provide a complete proof of
Theorem \ref{JNVN} in Section 4 below. When proving $(i)$, the polycyclic
case is rather simple while the non polycyclic case is more delicate and
relies on J. Groves' paper.

In \cite{BG}, we show that a similar phenomenon as in Corollary \ref{cor}
occurs for all finitely generated subgroups of the linear group GL$_{d}.$ It
would be interesting to study this also in the more general case of
elementary amenable groups.

\medskip

However, there is no uniform bound on $d_{\Gamma }^{+}$ itself in Theorem 
\ref{main}. As Osin pointed it out in \cite{Os2}, making use of a
construction of Grigorchuk and de la Harpe \cite{GH}, it is possible to
construct, for every $\varepsilon >0,$ an elementary amenable group $%
G_{\varepsilon }$ (and even solvable virtually metabelian, with polycyclic
and non polycyclic examples, as Bartholdi and Cornulier recently verified in 
\cite{BdC}) such that $0<S_{G_{\varepsilon }}<\varepsilon $.

The sharp contrast between polycyclic and non polycyclic groups is well
illustrated by the special case of metabelian groups. Indeed we have $%
d_{\Gamma}^{+} \leq 3$ for every non polycyclic metabelian group (this can
be derived both from Osin's method or ours, see Proposition \ref{nonpolyc})
while the following example shows that there is no upper bound for $%
d_{\Gamma}^{+}$ in the class of polycyclic metabelian groups.



\begin{theorem}
\label{ct}There exists a sequence $(G_{n})_{n\geq 1}$ of metabelian
polycyclic subgroups of $GL_{2}(\overline{\Bbb{Q}})$ such that $%
d_{G_{n}}^{+}\geq n$ for all $n.$
\end{theorem}

However this does not rule out the possibility that $S_{\Gamma}$ may be
bounded away from $0$ by a universal bound for all metabelian groups that
are not virtually nilpotent. Following Osin \cite{Os2} (remark after Theorem
2.4) let us formulate this as a question:

\begin{ques}
\label{ques}Is it true that there exists an $\varepsilon _{0}>0$ such that $%
S_{\Gamma }\geq \varepsilon _{0}$ for all non virtually nilpotent finitely
generated metabelian groups ?
\end{ques}

As any such group maps to the $2\times 2$ matrices by Theorem \ref{JNVN}
(ii), the question reduces to subgroups of $GL_{2}(\overline{\Bbb{Q}})$ and
even to the $2$-generated groups $\Gamma (x)$ defined in Section 7 below. It
is interesting to observe that a positive answer to this question would
imply the famous Lehmer conjecture from number theory, we explain this in
Section 7.

\bigskip

\textbf{Outline of the paper: }In Section 2, we explain how to build two
positively independent elements via the ping-pong argument. Section 3 is
devoted to the proof of Theorem \ref{main} in the particular case of
subgroups of the affine group. In Section 4, we provide a complete proof of
Theorem \ref{JNVN} along the lines of J. Groves' paper \cite{Gro}. In
Section 5, we complete the proof the Theorem \ref{main}. Section 6 deals
with elementary amenable groups and we then prove Theorem \ref{first}.
Finally in Section 7 we give examples of simple metabelian groups with
arbitrarily large $d_{\Gamma }^{+},$ proving Theorem \ref{ct}.

\bigskip

\textbf{Remark} We recently (December 2006) learned from J.S. Wilson 
that he knew of this approach via Theorem \ref{JNVN} to show uniform 
exponential growth for solvable groups.

\section{Ping-pong on the affine line}

Let $\Bbb{A}$ be the algebraic group of affine transformations of the line $%
\{x\mapsto ax+b\}$, that is, if $K$ is any field, 
\begin{equation*}
\Bbb{A}(K)=\left\{ \left( 
\begin{array}{ll}
a & b \\ 
0 & 1
\end{array}
\right) ,a\in K^{\times },b\in K\right\}
\end{equation*}
Elements from $\Bbb{A}(K)$ are of two types: they can either fix a point and
hence be homotheties around that point (when $a\neq 1$), or fix none and be
pure translations (when $a=1$). We are now going to give a simple sufficient
condition for two elements in $\Bbb{A}(k)$ to generate a free semigroup when
the field $k$ is a local field (i.e. $\Bbb{R}$, $\Bbb{C}$, a finite
extension of $\Bbb{Q}_{p}$ or a field of Laurent series over a finite
field). The construction of a free semigroup often relies on the ping-pong
principle. This principle can take many different guises, one of which is
illustrated by the following lemma:

\begin{lemma}
(Ping-pong) \label{pingpong}Let $k$ be an archimedean (resp.
non-archimedean) local field. Let $x$ and $y$ be two affine transformations
of the $k$-line such that $x$ fixes $p\in k$ and $y$ fixes $q\in k$. Assume
moreover that $x$ acts by multiplication by $\alpha $ around $p$ while $y$
acts by multiplication by $\beta $ around $q.$ Suppose that $|\alpha |$ and $%
|\beta |$ are $\leq \frac{1}{3}$ (resp. $<1$). Then $x$ and $y$ generate a
free semigroup.
\end{lemma}

\proof%
Let $d=|p-q|$ and let $B(p)$ and $B(q)$ be the open balls of radius $d/2$
(resp. $d$ in the non-archimedean case) centered at $p$ and at $q.$ By
definition they are disjoint and it follows from the assumption on $\alpha $
and $\beta $ that $x$ maps both of them into $B(p)$ while $y$ maps both of
them into $B(q).$ Now suppose that $w_{1}=x^{n_{1}}y^{m_{1}}\cdot ...\cdot
x^{n_{k}}y^{m_{k}}$ and $w_{2}=x^{n_{1}^{\prime }}y^{m_{1}^{\prime }}\cdot
...\cdot x^{n_{l}^{\prime }}y^{m_{l}^{\prime }}$ are two non trivial words
in $x$ and $y$ with non-negative powers. If they give rise to the same
element of $H$, then they must be equal as abstract words, i.e. $k=l$, $%
n_{i}=n_{i}^{\prime }$, $m_{i}=m_{i}^{\prime }$ for all $i$'s. Indeed
multiplying by a negative power of $x$ on the right hand side if necessary,
we may assume that $n_{1}=0$ while $n_{1}^{\prime }>0$, but then when acting
on the affine line, $w_{1}$ would send $q$ to a point in $B(q)$, while $%
w_{2} $ would send $q$ to a point in $B(p),$ a contradiction.%
\endproof%

To play this game of ping-pong, we will need to be able to choose a suitable
embedding of a finitely generated field $K$ into an appropriate local field.
This is done via the following easy and classical fact:

\begin{lemma}
\label{tits}\label{TitsLemma}(see J. Tits \cite{Tit}) Let $K$ be a finitely
generated field and $\alpha \in K$. If $\alpha ^{-1}$ is not an algebraic
integer (i.e. over $\Bbb{Z}$ if $char(K)=0$ or over $\Bbb{F}_{p}$ if $%
char(K)=p$), then there exists an embedding $\sigma :K\hookrightarrow k$
into a non-archimedean local field $k$ with absolute value $|\cdot |_{k}$
such that $|\sigma (\alpha )|_{k}<1.$ Moreover there exists a positive
number $\varepsilon =\varepsilon (K)>0$ such that if $\alpha \in K$ is an
algebraic unit and satisfies $|\log |\sigma (\alpha )|_{k}|<\varepsilon $
for every embedding $\sigma :K\hookrightarrow k$ into an archimedean local
field $k$, then $\alpha $ is a root of unity.
\end{lemma}

\proof%
Let $\beta =\alpha ^{-1}$ and let $K_{0}$ be the prime field of $K$ (i.e. $%
\Bbb{Q}$ if $char(K)=0$ and $\Bbb{F}_{p}$ if $char(K)=p$)$.$ Suppose first
that $\beta $ is transcendental over $K_{0}$. Transcendental elements are
dense in every local field, so if $k$ is a non-archimedean local field
containing $K_{0},$ then one can indeed find a transcendental element $%
\sigma (\beta )$ in $k$ with $|\sigma (\beta )|>1.$ This gives rise to an
embedding $\sigma :K_{0}(\beta )\rightarrow k$ which can always be extended
to the whole of $K$ up to changing $k$ into a finite extension if necessary (%
$K$ is finitely generated). Now if $K_{0}(\beta )$ is an algebraic extension
and $|\beta |\leq 1$ for every non-archimedean absolute value $|\cdot |$ on $%
K_{0}(\beta )$, then $\beta $ must be an algebraic integer.

Finally suppose $\alpha $ is an algebraic unit. Then its degree over $K_{0}$
is bounded by a constant depending only on $K,$ namely by $[K:K_{0}(\zeta
_{1},...,\zeta _{r})]<+\infty $ where $K_{0}(\zeta _{1},...,\zeta _{r})$ is
a purely transcendental extension over which $K$ is algebraic (by Noether's
normalization theorem). However, there are only finitely many algebraic
units of given degree all of whose conjugates are bounded. Moreover,
Kronecker's theorem implies that if $|\delta |\leq 1$ for all conjugates $%
\delta $ of $\alpha $, then $\alpha $ must be a root of unity. Hence the
result. 
\endproof%

\section{Proofs for subgroups of the affine group $\Bbb{A}=\{x\mapsto ax+b\}$%
}

First consider a general finitely generated metabelian group $\Gamma $ with
the exact sequence 
\begin{equation*}
1\rightarrow M\rightarrow \Gamma \overset{\pi }{\rightarrow }Q\rightarrow 1
\end{equation*}
where $M$ and $Q$ are abelian groups. The group $Q$ acts on $M$ by
conjugation. If we denote by $A$ the subring of $End(M)$ generated by $Q,$
then $M$ becomes an $A$-module. The following is classical:

\textbf{Claim 1:} $M$ is finitely generated as an $A$-module.

Indeed let $\{x_{1},...,x_{n}\}$ be generators of $\Gamma $. Since $Q$ is a
finitely generated abelian group, it is finitely presented $Q=\left\langle
y_{1},...,y_{n}|r_{1},...,r_{m}\right\rangle $ where $y_{i}=\pi (x_{i}).$ We
then verify that $M$ is generated as an $A$-module by the $%
r_{i}(x_{1},...,x_{n})$'s and the $x_{i}$'s that already belong to $M$. This
proves the claim.

We are now going to prove Theorem \ref{main} and Proposition \ref{nonpolyc} 
\textit{with the additional assumption} that the group $\Gamma $ is a
subgroup of $\Bbb{A}(K)$ for some field $K$. In the next section, we will
prove the general case by reducing to this one. We can take $K$ to be the
field generated by the matrix coefficients of the generators of $\Gamma $ so
that $K$ will then be finitely generated. Let $T$ be the subgroup of $\Bbb{A}%
(K)$ consisting of pure translations. We take $M=\Gamma \cap T$ and let $\pi 
$ be the canonical projection map from $\Bbb{A}(K)$ to $K^{\times }$ so that 
$Q$ is viewed as a multiplicative subgroup of $K^{\times }.$

\subsection{Proof of Theorem \ref{main} for subgroups of the affine line}

\label{affpf}

Suppose now that $\Gamma $ is not virtually nilpotent. Then clearly $Q$ does
not lie within the group of roots of unity of $K$, because otherwise $Q$,
being finitely generated, would be finite, and $\Gamma $ would be virtually
abelian. Now suppose that $Q$ lies in the subgroup of $K^{\times }$
consisting of algebraic units. Then $K$ must have characteristic zero and
according to Lemma \ref{tits}, there exists a positive $\varepsilon $
depending only on the field $K,$ hence only on the group $\Gamma $, such
that, if $\alpha \in Q$ is not a root of unity then there exists an
embedding of $K$ into an archimedean local field $k$ such that $|\log
|\alpha |_{k}|>\varepsilon .$ Let $n_{0}$ be an integer (depending only on $%
\Gamma $) such that $e^{-n_{0}\varepsilon }<\frac{1}{3}.$ Then let $\Gamma
_{0}=\pi ^{-1}(Q^{n_{0}}).$ It is a subgroup of finite index in $\Gamma $,
because $Q^{n_{0}}=\{b^{n_{0}},b\in Q\}$ is a subgroup of finite index in
the finitely generated abelian group $Q$. Let $\Gamma _{1}$ be a subgroup of
finite index in $\Gamma _{0}.$ We want to show that $d_{\Gamma _{1}}^{+}\leq
3.$ Let $\{z_{1},...,z_{n}\}$ be generators of $\Gamma _{1}.$ There must be
at least one $z_{i},$ say $z,$ such that $\pi (z)$ is not a root of unity.
By definition of $\Gamma _{0},$ $z$ is of the form $z_{0}^{n_{0}}a$ where $%
z_{0}\in \Gamma $ and $a\in M.$ Hence $\pi (z)=\pi (z_{0})^{n_{0}}$ and $\pi
(z_{0})$ is not a root of unity either. By Lemma \ref{tits}, there must
exist an embedding of $K$ into an archimedean local field $k$ such that $%
|\log |\pi (z_{0})|_{k}|>\varepsilon .$ Up to changing $z$ into $z^{-1}$ if
necessary, we may assume that $|\pi (z)|_{k}<\frac{1}{3}.$ The element $z$
is a non trivial homothety on the affine line of the local field $k$. It
fixes a point $p\in k.$ If all other $z_{i}$'s fix the same point $p$ then $%
\Gamma _{1}$ lies in the stabilizer of $p,$ an abelian subgroup. Hence at
least one of the $z_{i}$'s, call it $w,$ satisfies $wp\neq p.$ Then the two
affine transformations $x:=z$ and $y:=wzw^{-1}$ satisfy the hypothesis of
Lemma \ref{pingpong} and hence generate a free semigroup. We are done.

Suppose now that $Q$ does not lie inside the group of algebraic units of $K.$
Then we show that in fact $d_{\Gamma _{1}}^{+}\leq 3$ for every finite index
subgroup $\Gamma _{1}$ in $\Gamma $. Indeed, let $\{x_{1},...,x_{n}\}$ be
generators of $\Gamma _{1}$. Then at least one of the $x_{i}$'s, say $x$, is
such that $\pi (x)$ is not an algebraic unit. Indeed $\pi (\Gamma _{1})$ is
of finite index in $Q$ and if it were contained in the group of algebraic
units of $K^{\times }$ then so would $Q,$ because if $\alpha ^{n}$ is a unit
then $\alpha $ is a unit too. Up to changing $x$ into $x^{-1}$, we may
assume that $\pi (x)^{-1}$ is not an algebraic integer. Then, according to
Lemma \ref{tits}, there exists an embedding of $K$ into a non-archimedean
local field $k$ such that $|\pi (x)|_{k}<1.$ The element $x$ must be a non
trivial homothety on the affine line of $k$. It fixes a point $p\in k.$
Again, all other $x_{i}$'s cannot fix the same point $p.$ Hence at least one
of the $x_{i}$'s, call it $w,$ satisfies $wp\neq p.$ Then the two affine
transformations $x$ and $y:=wxw^{-1}$ satisfy the hypothesis of Lemma \ref
{pingpong}. We are done.

This finishes the proof of Theorem \ref{main} in the special case considered
in the present section.

\subsection{Polycyclic versus non polycyclic in $\Bbb{A}(K)$}

Here we prove the following result that will be useful in the proof of
Proposition \ref{nonpolyc}. 

\begin{lemma}
\label{zoo}Let $K$ be a field and $\Gamma $ be a finitely generated subgroup
of $\Bbb{A}(K)$. Let $Q=\pi (\Gamma )$ be the image of $\Gamma $ under the
canonical projection homomorphism $\pi $ from $\Bbb{A}(K)$ to $K^{\times }$.
Then $\Gamma $ is polycyclic (resp. virtually nilpotent) if and only if $Q$
is contained in the subgroup of algebraic units of $K^{\times }$ (resp. the
subgroup of roots of unity).
\end{lemma}

\proof%
Note that the $Q$ action on $M$ comes from the action of $\Gamma $ by
conjugation on the normal subgroup $M$, and when $Q$ is viewed as a subgroup
of $K^{\times }$ and $M$ as an additive subgroup of $K$, then this action
coincides with the action by multiplication. Also we may assume that $K$ is
finitely generated, since $\Gamma $ is so. As observed in Lemma \ref{tits},
the subfield of $K$ consisting of algebraic elements over the prime field $%
K_{0}$ is a finite extension of $K_{0}.$

Suppose $\Gamma $ is polycyclic and let $\alpha \in \pi (\Gamma )$. Since $%
\alpha $ is arbitrary, it is enough to show that $\alpha $ is an algebraic
integer. Up to conjugating $\Gamma $ inside $\Bbb{A}(K)$, we may assume that 
$\Gamma $ contains the element $\gamma =\left( 
\begin{array}{ll}
\alpha & 0 \\ 
0 & 1
\end{array}
\right) .$ Conjugation by $\gamma $ acts on $M$ by multiplication by $\alpha 
$. Since $\Gamma $ is polycyclic, $M$ must be a finitely generated abelian
group. By the Cayley-Hamilton theorem, $\alpha $ must satisfy a polynomial
equation with coefficients in $\Bbb{Z}$. Hence $\alpha $ is an algebraic
integer.

Conversely, from claim 1 above, there is a finite set $a_{1},...,a_{n}$ in $%
M $ such that $M=Qa_{1}+...+Qa_{n}.$ Let $\mathcal{O}_{K}$ be the ring of
integers (over $\Bbb{Z})$ of $K.$ Then we get that $M$ lies inside $\mathcal{%
O}_{K}a_{1}+...+\mathcal{O}_{K}a_{n}.$ However $\mathcal{O}_{K}$ is a
finitely generated $\Bbb{Z}$-module, hence a finitely generated additive
group. It follows that $M$ itself is finitely generated and that $\Gamma $
must be polycyclic contrary to the hypothesis made on $\Gamma $. We are
done. The corresponding statement for virtually nilpotent groups is also
straightforward and we omit the proof. 
\endproof%

\section{Just non virtually nilpotent groups}

In this section, we prove Theorem \ref{JNVN}. A finitely generated solvable
group is called just non virtually nilpotent (JNVN) if it is not virtually
nilpotent but all of its proper quotients are virtually nilpotent. In \cite
{Gro}, making use of ideas of P. Hall, J. Groves studied the structure of
metanilpotent just non polycyclic groups (see also \cite{RobWil}). As it
turns out, his methods apply with only minor changes to JNVN groups.

\subsection{JNVN groups are virtually metabelian}

We first make a few observations; we refer the reader to the textbook \cite
{Rob} for more details. Note that any subgroup and any quotient of a
virtually nilpotent group is again virtually nilpotent. Also a finite index
subgroup of a JNVN group is again JNVN. Further note that if $N_{1}$ and $%
N_{2}$ are two non trivial normal subgroup of a JNVN group $G$, then $%
N_{1}\cap N_{2}\neq \{1\}$ because $G/N_{1}\cap N_{2}$ embeds in the product 
$G/N_{1}\times G/N_{2}$ which is virtually nilpotent. Since finitely
generated virtually nilpotent groups are \textit{max-n}, i.e. have the
maximal property on normal subgroups (any increasing sequence of normal
subgroups stabilizes), so are the JNVN groups. Given a solvable group $G,$
we denote by $Fit(G)$ the Fitting subgroup of $G$, i.e. the subgroup
generated by all nilpotent normal subgroups of $G$. It is a basic fact that
the subgroup generated by two normal nilpotent subgroups is again a normal
nilpotent subgroup. Hence if $G$ is max-n then $Fit(G)$ is itself nilpotent.
When $G$ is JNVN, we have more:

\begin{lemma}
\label{Fitabelian}Let $G$ be a JNVN group. Then $Fit(G)$ is abelian.
\end{lemma}

\proof%
Let $N$ be a normal nilpotent subgroup of $G$. If $N^{\prime }\neq \{1\}$
then $N/N^{\prime }$ is finitely generated (like any subgroup of a finitely
generated virtually nilpotent group). Hence $N$ itself is finitely generated
(see \cite{Rob} 5.2.17.). Therefore $G$, being an extension of two
polycyclic groups is itself polycyclic. We now show that $G$ must be
virtually nilpotent, leading the desired contradiction. We can find a
subgroup $H$ containing $N^{\prime },$ with finite index in $G,$ such that $%
H/N^{\prime }$ is nilpotent. Then $H$ must act unipotently by conjugation on 
$N/N^{\prime },$ i.e. 
\begin{equation}
(h_{1}-1)\cdot ...\cdot (h_{n}-1)N/N^{\prime }=\{1\}  \label{unipotent}
\end{equation}
for some $n$ any any $h_{i}$'s in $H.$ As $N$ acts trivially on $N/N^{\prime
}$, (\ref{unipotent}) also holds for all $h_{i}$'s in $HN.$ Hence $%
HN/N^{\prime }$ is nilpotent. We can apply Hall's criterion for nilpotence
(see \cite{Rob} 5.2.10.) which says that $HN$ too must be nilpotent. But $HN$
has finite index in $G$ so we are done. 
\endproof%

Furthermore:

\begin{lemma}
$Fit(G)$ is either torsion free or is a torsion group of prime exponent $p$.
\end{lemma}

\proof%
Since $G$ has max-n, the torsion subgroup of $Fit(G)$ has finite exponent,
hence is a finite direct product of the $p$-torsion factors. But as we have
mentioned above no two non trivial normal subgroups of $G$ intersect
trivially. Hence there must be at most one $p$-torsion factor, say $T_{p}$.
Again $pFit(G)$ and $T_{p}$ intersect trivially, hence either $T_{p}$ is
trivial and $Fit(G)$ is torsion free, or $pFit(G)$ is trivial and $%
Fit(G)=T_{p}.$ 
\endproof%

It is a basic property of polycyclic groups due to Malcev that they contain
a finite index subgroup whose derived group is nilpotent (see \cite{Rob}
Proposition 15.1.6.). Thus if $G$ is JNVN and polycyclic, then it contains a
finite index subgroup $H$ with $H^{\prime }\leq Fit(H),$ which means by
Lemma \ref{Fitabelian} above that $H$ is metabelian. So part $(i)$ of
Theorem \ref{JNVN} is proved in this case.

We now assume that $G$ is JNVN and non polycyclic. Up to passing to a normal
subgroup of finite index containing $Fit(G),$ we may assume that $G/Fit(G)$
is nilpotent (since the Fitting subgroup is a characteristic nilpotent
subgroup, the two Fitting subgroups actually coincide). In this setting we
are going to show that $G$ is metabelian.

Let us denote by $A$ the commutative ring equal to $\Bbb{Z}$ when $Fit(G)$
is torsion free and equal to $\Bbb{F}_{p}[T]$ when $Fit(G)$ is of exponent $%
p.$ Let $K$ be its field of fractions, i.e. $\Bbb{Q}$ or $\Bbb{F}_{p}(T)$
respectively. We make $Fit(G)$ into an $A$-module by letting $T$ act via
conjugation by some non trivial element $z\in G$ whose projection modulo $%
Fit(G)$ is of infinite order in the center of $G/Fit(G)$ (there always is
such an element because $G/Fit(G)$ is an infinite finitely generated
nilpotent group). The following lemma is crucial:

\begin{lemma}
$Fit(G)$ is a torsion free $A$-module and $F_{K}:=Fit(G)\otimes _{A}K$ is
finite dimensional over $K$.
\end{lemma}

\proof%
We follow \cite[Lemma~2]{Gro}. We first show that $Fit(G)$ is torsion free
as an $A$-module. This is clear if $A=\Bbb{Z}$. When $A=\Bbb{F}_{p}[T]$, the
torsion elements in $Fit(G)$ form a subgroup that can be written as an
increasing union of the (normal) subgroups killed by the the polynomial $%
T^{n!}-1$ (any polynomial in $\Bbb{F}_{p}[T]$ divides $T^{n!}-1$ for $n$
large enough). As $G$ has max-n, it follows that $Fit(G)$ is annihilated by
some $T^{n!}-1.$ This means that $z^{n!}$ commutes with $Fit(G)$. Hence the
subgroup generated by $Fit(G)$ and $z^{n!}$ is normal in $G$ and abelian
thus contradicting the definition of $Fit(G).$

We now turn to the second half of the statement. It is a consequence of a
result of P. Hall (see \cite{Rob} 15.4.3.) saying that $Fit(G)$ contains a
free $A$-module $F_{0}$ such that $Fit(G)/F_{0}$ is a torsion $A$-module
with non trivial anihilator. Let $r\in A\backslash \{0\}$ in the anihilator
of $Fit(G)/F_{0}$ and let $q\in A$ be relatively prime to $r.$ Then clearly $%
qFit(G)\cap F_{0}=qF_{0}.$ Hence $F_{0}/qF_{0}$ embeds in $Fit(G)/qFit(G)$
hence is finitely generated (because $qFit(G)$ is a non trivial normal
subgroup of $G$). Since $F_{0}$ is a free $A$-module, $F_{0}$ too must be
finitely generated. Hence $F_{0}\otimes _{A}K=Fit(G)\otimes _{A}K$ is finite
dimensional over $K$. 
\endproof%

We are now ready to complete the proof of Theorem \ref{JNVN} $(i)$. The
group $G$ acts by conjugation via its quotient $G/Fit(G)$ on $Fit(G)$ and
this action extends to a $K$-linear action on $F_{K}:=Fit(G)\otimes _{A}K$.
On the other hand, for any nilpotent subgroup $N$ of GL$_{n}(K)$, its
derived group $N^{\prime }$ acts unipotently on $K^{n},$ i.e. $%
(g_{n}-1)\cdot ...\cdot (g_{1}-1)=0$ for any $g_{1},...,g_{n}$ in $N^{\prime
}.$ Therefore $G^{\prime }$ acts unipotently on $F_{K}$ (here it is crucial
that $\dim _{K}F_{K}<+\infty $). Hence $(g_{n}-1)\cdot ...\cdot (g_{1}-1)f=0$
in $F_{K}$ for every $g_{i}$'s in $G^{\prime }$ and every $f\in Fit(G).$ As $%
Fit(G)$ is a torsion free $A$-module, it embeds naturally in $F_{K}$ and the
previous equality actually holds in $Fit(G).$ Hence $%
[g_{n},[g_{n-1},[...[g_{1},f]...]]]$ is trivial in $Fit(G).$ We conclude
that $G^{\prime }$ is nilpotent and that $G^{\prime }\leq Fit(G),$ that is $%
G/Fit(G)$ is abelian, so $G$ is metabelian.

\subsection{Embeddings of metabelian groups into the affine group}

In this paragraph, we prove the second part of Theorem \ref{JNVN}. More
precisely, we show the following:

\begin{prop}
\label{metabelian}Let $H$ be a finitely generated metabelian group. Then
there is a field $K$ and a homomorphism $\rho :H\rightarrow \Bbb{A}(K)$ with
the following property. If $H$ is not virtually nilpotent (resp. not
polycyclic), then $\rho (H)$ too is not virtually nilpotent (resp. not
polycyclic).
\end{prop}

Note that it is easy to see that every finitely generated metabelian group
can be embedded in some $\Bbb{A}(R)$ where $R$ is some commutative ring
(Magnus embedding, see also \cite{Rem}). The point here is to find a
suitable quotient of $R$ which is a field and preserves non virtually
nilpotence. For the proof, we could again refer to further results in J.
Groves' paper \cite{Gro} after passing to a JNVN quotient. However, for the
reader's convenience, we provide another, more constructive proof.

\proof%
The finitely generated metabelian group $H$ comes with the exact sequence 
\begin{equation*}
1\rightarrow M\rightarrow H\rightarrow Q\rightarrow 1
\end{equation*}
where $M$ and $Q$ are abelian groups. The group $Q$ acts on $M$ by
conjugation. If we denote by $A$ the (commutative) subring of $End(M)$
generated by $Q,$ then $M$ becomes an $A$-module and we have the natural map 
$\alpha :H\rightarrow Aut(M)$ sending $h$ to the conjugation by $h$ on $M.$
In \textit{Claim 1} above we showed that $M$ is a finitely generated $A$%
-module. Next we make the following observation:

\begin{lemma}
Let $\{h_{i}\}_{1\leq i\leq m}$ be a finite generating set for $H$, and let $%
z_{i}=\alpha (h_{i})\in A$. Suppose that there are positive integers $n_{i}$
and $k_{i}$ such that $(z_{i}^{n_{i}}-1)^{k_{i}}=0$ in $A$ for all $i=1,...,m
$. Then $H$ is virtually nilpotent.
\end{lemma}

\proof%
Indeed, if this were true, there would be $n\geq 1$ and $k\geq 1$ such that $%
(z_{i}^{n}-1)^{k}=0$ for all $i=1,...,m$. Now let $Q_{0}$ be the subgroup of 
$Q$ generated by the $\pi (h_{i}^{n})$ for $i=1,...,m$. It has finite index
in $Q$ and its pull-back $H_{0}:=\pi ^{-1}(Q_{0})$ has finite index in $H$.
Clearly $H_{0}$ is the subgroup of $H$ generated by $M$ and the $h_{i}^{n}$%
's. Let $y_{i}=z_{i}^{n}.$ Then $[H_{0},H_{0}]\subset M$ while 
\begin{equation*}
\lbrack H_{0},M]=\left\langle (P_{1}-1)a|P_{1}\in \mathcal{P},a\in
M\right\rangle
\end{equation*}
where $\mathcal{P}$ is the set of all monomials in $y_{1},...,y_{m},$ i.e. $%
P_{1}=y_{1}^{k_{1}}\cdot ...\cdot y_{m}^{k_{m}}$ for some $k_{i}\geq 0.$
Similarly with $d$ commutators, 
\begin{equation*}
\lbrack H_{0},[H_{0}...[H_{0},M]...]=\left\langle (P_{d}-1)\cdot ...\cdot
(P_{1}-1)a|P_{1},...,P_{d}\in \mathcal{P},a\in M\right\rangle
\end{equation*}
Moreover for every $P\in \mathcal{P}$ there are elements $r_{i}\in A$ such
that $P-1=\sum r_{i}(y_{i}-1).$ Hence if $d>(k-1)m$ every product of the
form $(P_{d}-1)\cdot ...\cdot (P_{1}-1)$ can be written as a sum $\sum
r_{i}^{\prime }(y_{i}-1)^{k}$ for some $r_{i}^{\prime }\in A,$ hence is
always zero by the hypothesis on the $y_{i}$'s. It follows that $H_{0}$ is
nilpotent of order $d+1$ at most. Hence $H$ is virtually nilpotent. 
\endproof%

Similarly, it is (even more) straightforward to see that:

\begin{lemma}
Keeping the notation of the previous lemma, suppose that there are
polynomials $q_{i}\in \Bbb{Z}[X]$ with leading coefficient equal to $1$ such
that $q_{i}(z_{i})=0$ in $A$ for all $i=1,...,m.$ Then $H$ is polycyclic.
\end{lemma}

Let us go back to the proof of Proposition \ref{metabelian}. According to
the lemmas above, we can choose one of the $z_{i}$'s, call it $z$, such that 
$(z^{n}-1)^{k}\neq 0$ in $A$ for all positive integers $n$ and $k$. And if $%
H $ is not polycyclic we may even assume that $q(z)\neq 0$ in $A$ for all
monic polynomials $q\in \Bbb{Z}[X].$

Let $S$ be the subset of $A$ consisting of all products of factors of the
form $\phi (z)$, where $\phi \in \Bbb{Z}[X]$ runs through the collection of
all cyclotomic polynomials, i.e. $\{1,X-1,X+1,X^{2}+X+1,...\}.$ When $H$ is
not polycyclic we change the definition of $S$ and consider instead the set
of all $q(z)$ where $q\in \Bbb{Z}[X]$ is an arbitrary monic polynomial.
Clearly $S$ is a multiplicative subset of $A$ and, according to the choice
of $z$, $S$ does not contain $0$.

We can thus consider the localized ring $S^{-1}A.$ Then $S^{-1}M$ is a non
trivial finitely generated $S^{-1}A$-module. Let $I$ be the annihilator of $%
S^{-1}M,$ i.e. $I=\{r\in S^{-1}A,$ $rS^{-1}M=0\}$. Then $I$ is a proper
ideal if $S^{-1}A.$ It is therefore contained in a maximal ideal $P$ of $%
S^{-1}A.$ We can now set $K=S^{-1}A/P$ and $M_{0}=S^{-1}M/PS^{-1}M$.

$K$ is a field and $M_{0}$ a finite dimensional $K$-vector space. Moreover, $%
M_{0}$ is non trivial. Indeed, suppose that $S^{-1}M=PS^{-1}M.$ Then $%
(S^{-1}M)_{P}=P(S^{-1}M)_{P}$ as $(S^{-1}A)_{P}$-modules, where $%
(S^{-1}A)_{P}$ is the local ring associated to the prime ideal $P\subset
S^{-1}A.$ Moreover $(S^{-1}M)_{P}$ is a finitely generated $(S^{-1}A)_{P}$%
-module. Hence Nakayama's lemma applies and shows that $(S^{-1}N)_{P}=0.$
But this means that there must exists $r\in S^{-1}A$ with $r\notin P$ such
that $rS^{-1}M=0$ and contradicts the choice of $P$.

We are now in a position to define the desired linear representation. Let $%
\pi $ be either one of the canonical maps $A\rightarrow K$ or $M\rightarrow
M_{0}.$ Let $h_{0}$ be an element of $H$ such that $\alpha (h_{0})=z$, and
define the map $c:H\rightarrow M_{0}$ by 
\begin{equation*}
c(h)=\pi ([h_{0},h])
\end{equation*}
It is straightforward to check that for all $h_{1}$ and $h_{2}$ in $H$ 
\begin{equation*}
c(h_{1}h_{2})=\beta (h_{1})c(h_{2})+c(h_{1})
\end{equation*}
where $\beta =\pi \circ \alpha .$

As $z-1$ is invertible in $S^{-1}A,$ we see that $c(M)=M_{0}.$ Let $%
z_{0}:=\pi (z).$ Then, according to the choice to $z$, $z_{0}\in K^{\times }$
is not a root of unity, and in the case when $H$ is not polycyclic $z_{0}$
is even not an algebraic integer in $K$.

Let $F_{0}$ be a hyperplane in $M_{0}$ and let $\pi _{0}:$ $M_{0}\rightarrow
M_{0}/F_{0}\simeq K$. Setting $\overline{c}=c\circ \pi _{0},$ we obtain the
desired representation of $H$ into the affine group $\Bbb{A}(K)$ defined by 
\begin{eqnarray*}
\rho &:&H\rightarrow \Bbb{A}(K) \\
h &\mapsto &\left( 
\begin{array}{ll}
\beta (h) & \overline{c}(h) \\ 
0 & 1
\end{array}
\right)
\end{eqnarray*}
And this homomorphism sends $h_{0}$ to $\left( 
\begin{array}{ll}
z_{0} & 0 \\ 
0 & 1
\end{array}
\right) $ where $z_{0}\in K^{\times }$ is not a root of unity (resp. not an
algebraic integer when $H$ is non polycyclic), and $\rho (M)$ equals $\left( 
\begin{array}{ll}
1 & K \\ 
0 & 1
\end{array}
\right) $. By Lemma \ref{zoo}, any such subgroup of $\Bbb{A}(K)$ is not
virtually nilpotent (resp. non polycyclic). We are done. 
\endproof%

\section{Proof of Theorem \ref{main}}

To prove Theorem \ref{main}, it is enough to show that some quotient of $%
\Gamma $ has the same property. Combining the proof of these results we gave
in Section 3 for subgroups of $\Bbb{A}(K)$ with Theorem \ref{JNVN}, we see
that we are done if we make use of the following easy and standard fact:

\textbf{Claim 2:} Every finitely generated non virtually nilpotent group has
a just non virtually nilpotent quotient.

Indeed, let $G$ be such a group and let $\mathcal{N}$ be the set of all
normal subgroups $N$ of $G,$ such that $G/N$ is not virtually nilpotent.
Suppose $N_{1}\subset N_{2}\subset ...\subset N_{i}\subset ...$ is an
increasing chain of subgroups from $\mathcal{N}$. And let $N$ be the union
of all $N_{i}$'s. Then $N$ is indeed a normal subgroup of $G.$ Now if $G/N$
were virtually nilpotent, there would exist a subgroup of finite index $%
G_{0} $ in $G,$ containing $N,$ such that $G_{0}/N$ is nilpotent. Like any
finitely generated nilpotent group, $G_{0}/N$ has a finite presentation $%
\left\langle x_{1},...,x_{n}|r_{1},...,r_{m}\right\rangle .$ The finitely
many relations $r_{i}$'s belong to one of the $N_{i}$'s, say $N_{i_{0}}.$
Hence $G_{0}/N_{i_{0}}$ appears as a quotient of $G_{0}/N,$ hence is
nilpotent, contradicting the assumption that $G/N_{i_{0}}$ is not virtually
nilpotent. It follows that we can apply Zorn's lemma and obtain a maximal
element $N$ in $\mathcal{N}.$ Then clearly $G/N$ is not virtually nilpotent,
while any proper quotient of it is. qed.

\section{Elementary amenable groups}

In this section, we explain why Osin's work actually implies Theorem \ref
{first} and also Proposition \ref{Osin} below.

We first discuss the following lemma:


\begin{lemma}
\label{pptable}Let $\Gamma $ be a group generated by two elements $x$ and $y$
such that the normal subgroup $M=\left\langle x^{n}yx^{-n},n\in \Bbb{Z}%
\right\rangle $ is not finitely generated. Then the elements $x$ and $%
yxy^{-1}$ generate a free semigroup.
\end{lemma}

It is easy to prove this lemma directly by showing that equality between two
different positive words in $x$ and $yxy^{-1}$ would force $M$ to be
generated by finitely many $x^{n}yx^{-n}$, we leave this exercise to the
reader. A. Navas pointed out to me that this lemma/exercise is also stated
(with its solution) in \cite{Rhem}. Also, as the referee pointed out to me,
a careful analysis of Milnor's argument in \cite{Mil} shows that it is in a
way already contained there. However it is also possible to give a different
(and more complicated!) proof using yet another ping-pong argument! And we
now explain this idea.

\proof%
Let $M_{n}=\left\langle x^{k}yx^{-k},n\leq k\right\rangle $ and $%
M_{n}^{-}:=\left\langle x^{k}yx^{-k},n\geq k\right\rangle $. Suppose $y\in
M_{1}\cap M_{-1}^{-}.$ Then $M$ would be finitely generated. Indeed there
exists $N\in \Bbb{N}$ such that $y\in \left\langle x^{k}yx^{-k},1\leq k\leq
N\right\rangle $ . Hence all $x^{k}yx^{-k}$ for negative $k$ also belong to $%
\left\langle x^{i}yx^{-i},1\leq i\leq N\right\rangle .$ Similarly all the $%
x^{k}yx^{-k}$ for positive $k$ would belong to $\left\langle
x^{i}yx^{-i},-N^{\prime }\leq i\leq -1\right\rangle $ for some $N_{1}\in 
\Bbb{N}$. Therefore up to changing $x$ into $x^{-1}$ we may assume that $%
y\notin M_{\infty }=\cap _{n\in \Bbb{Z}}M_{n}.$ Now $\Gamma $ acts on $%
M/M_{\infty }$ via the action $\gamma \cdot aM_{\infty
}=x^{n}bax^{-n}M_{\infty }$ where $\gamma =x^{n}b$ for some $n\in \Bbb{Z}$
and $b\in M.$

We define a valuation $v$ on $M$ by $v(a)=\sup \{n,$ $a\in M_{n}\}$ --
clearly $v(ab)\geq \min \{v(a),v(b)\}$ -- and an ultrametric distance $%
d(x,y)=e^{-v(a^{-1}b)}.$ This is a well defined distance on the quotient
space $M/M_{\infty }.$ Let $X$ be the completion of $M/M_{\infty }$ for the
distance $d$. Then the $\Gamma $-action extends to an action on $X.$
Moreover $M$ acts by isometries while $x$ contracts distances by a factor $%
e^{-1}$ and $x^{-1}$ is $e$-Lipschitz. It is easy to check that the elements
of $\Gamma $ may be of three possible kinds: if $\gamma \notin M$ then $%
\gamma $ has a unique fixed point on $X$ and acts by contraction/dilatation
by a facteur $e^{-1}$ around it (the sequence $x_{n+1}=\gamma \cdot x_{n}$
is a Cauchy sequence), if $\gamma \in M$ then either $\gamma $ has no fixed
point and acts like a translation or it fixes a point and belongs to a
conjugate of $M_{\infty }.$ Clearly $x$ is of the first kind, fixes the
coset $M_{\infty },$ while $y$ doesn't.

Just like in the non-archimedean case of Lemma \ref{pingpong}, we see that $%
x $ and $yxy^{-1}$ are two contractions with different fixed points and play
ping-pong on $X.$ 
\endproof%

In \cite{Os2} Osin also shows a strong uniformity statement about the growth
of infinitely--generated--by--nilpotent groups. His argument actually gives
the following bound for $d_{\Gamma}^{+}$.

\begin{proposition}
\label{Osin}Let $\Gamma $ be a finitely generated group given by an
extension 
\begin{equation*}
1\rightarrow M\rightarrow \Gamma \rightarrow N\rightarrow 1
\end{equation*}
where $N$ is $r$-step nilpotent and $M$ is not finitely generated. Then $%
d_{\Gamma }^{+}\leq 3^{2r+1}.$
\end{proposition}

\begin{proof}
 In \cite{Os2}, Osin proves using commutator calculus and an intricate
induction the remarkable fact that if $S$ is a symmetric generating set for $%
\Gamma $ and all subgroups of the form $\left\langle x^{n}yx^{-n},n\in \Bbb{Z%
}\right\rangle $ where $x\in S^{(r)}$ and $y\in S^{(2r)}$ are finitely
generated ($S^{(n)}$ denotes the finite set of all commutators of length at
most $n$ in the elements of $S$) then $M$ too is finitely generated. See 
\cite[Section 5 and Lemma 6.4.]{Os2} The set $S^{(r)}$ lies in the ball of
radius $3^{r}$ (rough estimate) for the word metric induced by $S$. Hence we
can apply Lemma \ref{pptable} to some choice of $x$ and $y.$ Then $x$ and $
yxy^{-1}$ will lie in the ball of radius $3^{2r+1}.$%
\end{proof}

We now complete the proof of Theorem \ref{first}. 
\begin{proof}[Proof of Theorem \ref{first}]
 Again as Osin observes in \cite[Proposition~3.1.]{Os2}, if $\Gamma $ is a
finitely generated elementary amenable group which is just non virtually
nilpotent, then either $\Gamma $ is virtually polycyclic or $\Gamma $ has a
non trivial normal subgroup $M$ which is not finitely generated. We can pass
to a subgroup of finite index $\Gamma _{0}$, since $d_{\Gamma }^{+}\leq
(2[\Gamma :\Gamma _{0}]+1)\cdot d_{\Gamma _{0}}^{+}.$ The polycyclic case
was already treated in Theorem \ref{main} and the other case follows
directly from Proposition \ref{Osin} above. 
\end{proof}

\section{Metabelian groups, a counter-example, and the Lehmer conjecture}
\label{ctrex}\label{cnterex}

As mentionned in the introduction we have the following:

\begin{proposition}
\label{nonpolyc}Let $\Gamma $ be a finitely generated metabelian group which
is not polycyclic. Then $d_{\Gamma }^{+}\leq 3.$
\end{proposition}

\begin{proof}
By \textit{Claim 2} in Section 5 and Theorem \ref{JNVN} (ii) we can assume that $\Gamma$
is a subgroup of $\Bbb{A}(K)$ for some field $K$. By Lemma \ref{zoo}, the quotient group 
$Q$ does not lie in the group of algebraic units of $K$. Hence we are in the situation
 of \ref{affpf} at the end of the proof of Theorem \ref{main} for subgroups of $\Bbb{A}(K)$,
 hence $d_{\Gamma}^{+}\leq 3$.
\end{proof}

Next we prove Theorem \ref{ct} from the introduction, namely:

\begin{theorem}
For every integer $n\geq 1$, there exists a $2$-generated polycyclic
subgroup $G_{n}$ of the affine group $\Bbb{A}(\overline{\Bbb{Q}})$ such that 
$d_{G_{n}}^{+}\geq n$ and $G_{n}$ is not virtually nilpotent.
\end{theorem}

\proof%
For $x\in \Bbb{C}$, let $\Gamma (x)$ be the subgroup of $\Bbb{A}(\Bbb{C})$
generated by the matrices $\left( 
\begin{array}{ll}
x & 0 \\ 
0 & 1
\end{array}
\right) $ and $\left( 
\begin{array}{ll}
1 & 1 \\ 
0 & 1
\end{array}
\right) $ that is the set of matrices of the form $\left( 
\begin{array}{ll}
x^{n} & P(x) \\ 
0 & 1
\end{array}
\right) $ where $n\in \Bbb{Z}$ and $P\in \Bbb{Z}[X,\frac{1}{X}]$. Note that $%
\Gamma (x)$ is also generated by $A(x)=\left( 
\begin{array}{ll}
x & 0 \\ 
0 & 1
\end{array}
\right) $ and $B(x)=\left( 
\begin{array}{ll}
x & 1 \\ 
0 & 1
\end{array}
\right) .$ We are going to exhibit a sequence of points $x_{n}\in \Bbb{C}$
such that $G_{n}:=\Gamma (x_{n})$ satisfies the requirement of the theorem,
that is more precisely $d^{+}(\Sigma _{n})\geq n$ where $\Sigma
_{n}=\{A(x_{n}),B(x_{n})\}$.

First, we make the following observation. Let $a,b,a^{\prime },b^{\prime }$
be non trivial homotheties in $\Bbb{A}(\Bbb{C})$ with respective dilation
ratio $x,y,x^{\prime },y^{\prime }\in \Bbb{C}^{\times }.$ Assume that $a$
and $b$ have distinct fixed points and that so do $a^{\prime }$ and $%
b^{\prime }.$ Assume further that $x=x^{\prime }$ and $y=y^{\prime }.$ Then
the pair $\{a,b\}$ generates a free semigroup if and only if the pair $%
\{a^{\prime },b^{\prime }\}$ generates a free semigroup. Indeed we check
easily that the two pairs are conjugate by a single element $\gamma \in \Bbb{%
A}(\Bbb{C}),$ i.e. $a^{\prime }=\gamma a\gamma ^{-1}$ and $b^{\prime
}=\gamma b\gamma ^{-1}.$

Second, we note that if $a$ is a homothety with dilation ratio $x\in \Bbb{C}%
^{\times }$ and $b$ is a translation then the pair $\{a,b\}$ never generates
a free semigroup unless $x$ is transcendental. Indeed, suppose $x$ is
algebraic over $\Bbb{Q}$, i.e. it is a root of a polynomial $\pi(X)$ of
degree $d$ with coefficients in $\Bbb{Z}$, then it is straightforward to
check that the following non trivial relation is satisfied $%
b^{b_{0}}ab^{b_{1}}a...b^{b_{d}}=b^{a_{0}}ab^{a_{1}}a...b^{a_{d}}$, where $%
P(X)=a_0+...+a_dX^d$ and $Q(X)=b_0+...+b_dX^d$ are polynomials of degree $d$
with non-negative integer coefficients such that $\pi=P-Q$.

From this it follows that whether or not two elements from $\Gamma (x)$
generate a free semigroup is a property that depends only on the respective
dilation ratios of the two elements. And since in $\Gamma (x)$ all elements
have a dilation ratio of the form $x^{n}$ for some integer $n\in \Bbb{Z},$
we may define the subset $\mathcal{NF}(x)$ of $\Bbb{Z}^{2}$ to be the set of
all couples $(n,m)\in \Bbb{Z}^{2}$ such that pairs of elements of $\Gamma
(x) $ of ratio $x^{n}$ and $x^{m}$ respectively do not generate a free
semigroup. The goal is now to find $x_{n}\in \Bbb{C}$ such that $\mathcal{NF}%
(x_{n})$ contains all couples $(p,q)$ with $|p|,|q|\leq n$.

Note that the matrices $A(x)$ and $B(x)$ satisfy $A(x)^{4}=B(x)^{2}A(x)B(x)$
if and only if $x^{3}+x+1=0.$ Now let $x_{n}$ be a root of the equation $%
X^{3n!}+X^{n!}+1=0$. We can assume that $|x_{n}|<1$. Hence $%
A(x_{n}^{n!})^{4}=B(x_{n}^{n!})^{2}A(x_{n}^{n!})B(x_{n}^{n!})$. Note that
for every $p\in \Bbb{Z}$, the element $B(x_{n})^{p}$ does not fix $0$
because $x_{n}$ is not a root of unity. It follows that $B(x_{n}^{p})$ is
conjugate to $B(x_{n})^{p}$ by a element of $\Bbb{A}(\Bbb{C})$ that fixes $0$
(i.e. commutes with $A(x_{n})$). Hence for all integers $p$ and $q$ with $%
|p|,|q|\leq n$ we have $A(x_{n}^{q})^{\frac{4n!}{q}}=B(x_{n}^{p})^{\frac{2n!%
}{p}}A(x_{n}^{q})^{\frac{n!}{q}}B(x_{n}^{p})^{\frac{n!}{p}}.$ If follows
that $A(x_{n}^{q})$ and $B(x_{n}^{p})$ do not generate a free semigroup
(although they do not commute), therefore $(p,q)\in \mathcal{NF}(x_{n}).$ We
are done.%
\endproof%

Let us now come back to Question \ref{ques} from the introduction. Let us
show that a positive answer to the question would imply the Lehmer
conjecture from Number Theory. Recall that if $\pi \in \Bbb{Z}[X]$ is a
monic polynomial, and $\pi =\prod_{1\leq i\leq d}(X-\lambda _{i})$ its
factorization over $\Bbb{C}$, the Mahler measure of $\pi $ is the number 
\begin{equation*}
m(\pi )=\prod_{1\leq i\leq d}\max \{1,|\lambda _{i}|\}
\end{equation*}
A classical theorem of Kronecker says that $m(\pi )=1$ if and only if all
roots $\lambda _{i}$'s are roots of unity. The Lehmer conjecture states that
there exists a universal $\varepsilon _{0}>0$ such that if $m(\pi )>1$ then $%
m(\pi )>1+\varepsilon _{0}.$ Clearly, the conjecture reduces to the case
when $\pi $ is irreducible over $\Bbb{Q}$ and $\pi (0)=\pm 1.$ Let $x$ be a
root of such a $\pi $ and let us consider the group $\Gamma (x).$ The
following claim is what we are aiming for:

\textbf{Claim:}\ $log(m(\pi ))\geq S_{\Gamma (x)}$

Let $\Sigma =\{1,A(x)^{\pm 1},B(x)^{\pm 1}\}.$ Any element $w$ of $\Sigma
^{n}$ is of the form $\left( 
\begin{array}{ll}
x^{k} & P(x) \\ 
0 & 1
\end{array}
\right) $ where $|k|\leq n,$ $d^{\circ }P\leq n,$ and $\left\| P\right\|
\leq n$, where we have set $\left\| P\right\| =\sum |a_{i}|$ if $P=\sum
a_{i}X^{i}.$ Therefore $\#\Sigma ^{n}\leq (2n+1)\cdot \#\{P(x),d^{\circ
}P\leq n,$ and $\left\| P\right\| \leq n\}.$ Let $P=\pi Q+R_{P}$ the
Euclidean division of $P$ by $\pi .$ Let us give an upper bound on the
number of possible remainders $R_{P}$ for $P\in \Bbb{Z}[X]$ with $d^{\circ
}P\leq n,$ and $\left\| P\right\| \leq n.$ For $k\geq 0$, let $%
Y_{k}=R_{X^{k}}.$ Then it is clear that the coefficients of $Y_{k}$ in the
basis $1,X,...,X^{d-1}$ satisfy a linear recurrence relation, i.e. $%
Y_{k}=M^{k}Y_{0}$ where $M$ is the companion matrix of $\pi .$ Let $%
(v_{1},...,v_{d})$ be a basis $\Bbb{C}_{d-1}[X]$ diagonalizing $M$ ($M$ has
distinct eigenvalues, since $\pi $ is irreducible). If $Y_{0}=\Sigma _{1\leq
i\leq d}\alpha _{i}v_{i}$ is the expression of $Y_{0}$ in this basis, then $%
Y_{k}=\Sigma _{1\leq i\leq d}\alpha _{i}\lambda _{i}^{k}v_{i}$ for every $%
k\geq 1.$ Let $|\alpha |=\max_{1\leq i\leq d}\{|\alpha _{i}|\}$ and $%
|\lambda |=\max_{1\leq i\leq d}\{|\lambda _{i}|^{d-1}\}.$

Let $B_{n}=\left\{ \sum_{1\leq i\leq d}x_{i}v_{i},|x_{i}|\leq n|\alpha
|\cdot \max \{1,|\lambda _{i}|\}^{n}\text{ for each }i\right\} $ and $%
a_{n}(i)=n|\alpha |\cdot \max \{1,|\lambda _{i}|\}^{n}+\frac{d}{2}|\alpha
|\cdot |\lambda |$ and $B_{n}^{\prime }=\left\{ \sum_{1\leq i\leq
d}x_{i}v_{i}\in \Bbb{C}[X],|x_{i}|\leq a_{n}(i)\text{ for each }i\right\} .$
For each $P\in \Bbb{Z}[X]$ with $d^{\circ }P\leq n,$ and $\left\| P\right\|
\leq n$ we have $R_{P}\in B_{n}\cap \Bbb{Z}[X],$ and if $x\in B_{n}\cap \Bbb{%
Z}[X],$ then $x+\varepsilon \in B_{n}^{\prime }$ for every $\varepsilon
=\sum_{1\leq i\leq d}\varepsilon _{i}X^{i-1}$ with $|\varepsilon _{i}|\leq 
\frac{1}{2}.$ Therefore $\#\{R_{P},P\in \Bbb{Z}[X],d^{\circ }P\leq n,\left\|
P\right\| \leq n\}\leq vol(B_{n}^{\prime }\cap V)$, where $vol$ is the
standard Lebesgue measure in $V=\Bbb{R}_{d-1}[X]$ in the basis $%
Y_{0}=1,Y_{1}=X,...,Y_{d-1}=X^{d-1}.$ But if $A_{n}$ is the endomorphism of $%
\Bbb{C}_{d-1}[X]$ such that $A_{n}(v_{i})=v_{i}\cdot a_{n}(i)$ then $%
vol(B_{n}^{\prime }\cap V)=vol(A_{n}(B^{\prime }\cap
O_{n}V))=\det_{(Y_{0},...,Y_{d-1})}(A_{n})vol(B^{\prime }\cap O_{n}V),$
where $B^{\prime }=\left\{ \sum_{1\leq i\leq d}x_{i}v_{i}\in \Bbb{C}%
[X],|x_{i}|\leq 1\right\} $ and $O_{n}$ is an orthogonal transformation such
that $O_{n}V=A_{n}^{-1}V$. Hence 
\begin{equation*}
vol(B_{n}^{\prime }\cap V)=\left( \prod_{1\leq i\leq d}a_{n}(i)\right) \cdot
\det_{(v_{1},...,v_{d})}(Y_{0},...,Y_{d-1})\cdot vol(B^{\prime }\cap O_{n}V)
\end{equation*}
But $vol(B^{\prime }\cap O_{n}V)$ converges to a non zero limit, as $O_{n}V$
converges in the grassmannian. Hence 
\begin{equation*}
\lim_{n\rightarrow +\infty }\left( vol(B_{n}^{\prime }\cap V)\right) ^{\frac{%
1}{n}}=\lim_{n\rightarrow +\infty }\left( \prod_{1\leq i\leq
d}a_{n}(i)\right) ^{\frac{1}{n}}=m(\pi )
\end{equation*}
Finally $\overline{\lim_{n\rightarrow +\infty }}\left( \#\Sigma ^{n}\right)
^{\frac{1}{n}}\leq m(\pi ),$ which is what we wanted. 
\endproof%

\begin{Ack}
I would like to thank Yves Guivarc'h for drawing my attention to his early
work \cite{Gui0}, Yves de Cornulier for pointing out to me the references 
\cite{Os3}, \cite{Rem} and \cite{RobWil}, Pierre de la Harpe and Andres
Navas for their remarks on an earlier version of the paper.
\end{Ack}


\begin{thebibliography}{99}
\bibitem{Alp}  R. Alperin, \textit{Uniform exponential growth of polycyclic
groups}, Geom. Dedicata \textbf{92} (2002), p. 105--113.

\bibitem{BdC}  L. Bartholdi, Y. de Cornulier, \textit{Infinite groups with
large balls of torsion elements and small entropy}, to appear in Archiv der
Mathematik.

\bibitem{BG}  E. Breuillard, T. Gelander, \textit{A uniform Tits alternative}%
, preprint.

\bibitem{Chou}  C. Chou, \textit{Elementary amenable groups}, Illinois J.
Math. \textbf{24} (1980), no. 3, 396--407.

\bibitem{EMO}  A. Eskin, S. Mozes, H. Oh, \textit{On uniform exponential
growth for linear groups}, Invent. Math. \textbf{160} (2005), no. 1, 1--30.

\bibitem{GH}  R. Grigorchuk, P. de la Harpe, \textit{Limit behaviour of
exponential growth rates for finitely generated groups}, Essays on geometry
and related topics, Vol. 1, 2, 351--370, Monogr. Enseign. Math., \textbf{38}%
, Enseignement Math., (2001).

\bibitem{Gro}  J. Groves, \textit{Soluble groups with every proper quotient
polycyclic}, Illinois J. Math. \textbf{22} (1978), no. 1, 90--95.

\bibitem{Grom}  M. Gromov, \textit{Hyperbolic groups. Essays in group theory}%
, 75--263, Math. Sci. Res. Inst. Publ., \textbf{8}, Springer, (1987).

\bibitem{Gui0}  Y. Guivarc'h, \textit{G\'{e}n\'{e}rateurs des groupes r\'{e}%
solubles}, Publications des S\'{e}minaires de Math\'{e}matiques de
l'Universit\'{e} de Rennes (1968), Ann\'{e}e 1967--1968, Exp. No. \textbf{1}%
, 17pp.

\bibitem{Har}  P. de la Harpe, \textit{Topics in Geometric Group Theory},
Chicago University Press, (2001).

\bibitem{Rhem}  P. Lomgobardi, M. Maj, A.H. Rhemtulla, \textit{Groups with
no Free Subsemigroups}, Trans. Amer. Math. Soc., Vol \textbf{347}, No 4,
(1995) p. 1419--1427.

\bibitem{Mil}  J. Milnor, \textit{Growth of finitely generated solvable
groups}, J. Diff. Geometry \textbf{2} (1968) p. 447-449.

\bibitem{Os}  D. Osin, \textit{The entropy of solvable groups}, Erg. Theory.
Dyn. Sys. \textbf{23}, no. 3, (2003) p. 907--918.

\bibitem{Os2}  D. Osin, \textit{Algebraic entropy of elementary amenable
groups}, Geometriae Dedicata 107, (2004) p. 133--151.

\bibitem{Os3}  D. Osin, \textit{Uniform non-amenability of free Burnside
groups}, preprint arXiv math.GR/0404073.

\bibitem{Rem}  V. N. Remeslennikov, \textit{Representation of finitely
generated metabelian groups by matrices}, Algebra i Logika \textbf{8} (1969)
p. 72--75.

\bibitem{Rob}  D. Robinson, \textit{A course in the theory of groups},
Springer Verlag, GTM.

\bibitem{RobWil}  D. Robinson, J. Wilson, \textit{Soluble groups with many
polycyclic quotients}, Proc. London Math. Soc. (3) \textbf{48} (1984), no.
2, p. 193--229.

\bibitem{Ros}  J. Rosenblatt, \textit{Invariant measures and growth
conditions}, Trans. Amer. Math. Soc. \textbf{193} (1974), p. 33--53.

\bibitem{Wil}  J. Wilson, \textit{On exponential growth and uniformly
exponential growth for groups}, Invent. Math. \textbf{155} (2004), no. 2, p.
287--303.

\bibitem{W}  J. Wolf,\textit{\ Growth of finitely generated solvable groups
and curvature of Riemanniann manifolds}, J. Differential Geometry, \textbf{2}
(1968) p. 421--446.

\bibitem{Tit}  J. Tits, \textit{Free subgroups in linear groups}, J. Algebra 
\textbf{20} (1972) p. 250-270.
\end{thebibliography}
\end{document}